\documentclass[12pt]{article}
\usepackage{a4wide}
\usepackage{amssymb}
\usepackage{amsmath}
\usepackage[russian]{babel}

  \chardef\No=242
\newtheorem{theorem}{Theorem}{}
\newtheorem{definition}{Definition}{}
{}
\newtheorem{lemma}{Lemma}{}
\newtheorem{remark}{Remark}{}

\begin{document}

\begin{center}
{\bf On a weighted variable spaces $L_{p(x),\,\omega}$ for $0< p(x)< 1$\\
and weighted Hardy inequality}
\end{center}
\begin{center}
ROVSHAN A.BANDALIEV
\end{center}

{\small  ABSTRACT. In this paper a weighted variable exponent Lebesgue spaces
$L_{p(x),\,\omega}$ for $0< p(x)< 1$ is investigated. We show that this spaces is
a quasi-Banach spaces. Note that embedding theorem between weight variable Lebesgue spaces is proved.
In particular, we show that $L_{p(x),\,\omega}(\Omega)$ for $0< p(x)< 1$ isn't locally convex.
Also, in this paper a some two-weight estimates for Hardy operator are proved.}

\vspace{2mm}


\vspace{4mm}
\begin{center}
1. {\bf  Introduction.}
\end{center}

It is well known that the variable exponent Lebesgue space $L_{p(x)}$ for $p(x)\ge 1$ appeared
in the literature for the first time already in [13]. Further development of this
theory was connected with the theory of modular function spaces. Somewhat later, a
more explicit version of these spaces, namely modular function spaces, were investigated
by many mathe-maticians (see [12]). The next step in the investigation of variable exponent
spaces was given in [16] and in [8]. But the variable exponent Lebesgue space for $0< p(x)< 1$
very less studied. Note that the space $L_{p(x)}$ for $0< p(x)< 1$ isn't modular function spaces.
The study of these spaces has been stimulated by problems of elasticity,
fluid dynamics, calculus of variations and differential equations with non-standard growth
conditions$\,$ (see [14], [17],[18]). For detailed information about variable exponent Lebesgue
space $L_{p(x)}$ for $p(x)\ge 1$ we refer to [7].

Let $R^{n}$ be the $n$-dimensional Euclidean space of points $x=\left(x_{1},..., x_{n}\right)$
and $\Omega$ be a Lebesgue measurable subset in $R^n$ and $\displaystyle{|x|= \left(\sum\limits_{i= 1}^n
x_i^2\right)^{1/2}}.$ Suppose that $p$ is a Lebesgue measurable function on $\Omega$ such that
$0< \underline p\le p(x)\le \overline p< 1,$ $\underline p= ess\,\inf_{x\in \Omega} p(x),$
$\overline p= ess\,\sup_{x\in \Omega} p(x),$ and $\omega$ is a weight function on $\Omega,$ i.e.
$\omega$ is non-negative, almost everywhere (a.e.) positive function on $\Omega.$ The Lebesgue measure
of a set $\Omega$ will be denoted by $|\Omega|.$ It is well known that $\displaystyle{|B(0,1)|= \frac
{\pi^{\frac n2}}{\Gamma\left(\frac n2+ 1\right)}},$ where $B(0,1)= \left\{x:\,x\in R^n;\;|x|< 1\right\}.$
Further, in this paper all sets and functions are supposed Lebesgue measurable.

\vspace{3mm}
\begin{center}
2. {\bf Preliminaries}
\end{center}

\begin{definition}
By $L_{p(x),\,\omega}(\Omega)$ we denote the set of measurable
functions $f$ on $\Omega$ such that
$$
I_{p,\,\omega}(f)= \int\limits_{\Omega} \left(|f(x)|\,\omega(x)\right)^{p(x)}\,dx< \infty.
$$
Note that the expression
$$
\|f\|_{L_{p(x),\,\omega}(\Omega)}= \|f\|_{p,\,\omega,\,\Omega}=
\inf\left\{\lambda> 0:\;\;\int\limits_{\Omega} \left(\frac{|f(x)|\,\omega(x)}
{\lambda}\right)^{p(x)} \,dx\le 1\right\} \eqno(2.1)
$$
defines a quasi-Banach spaces.
\end{definition}

We note some main properties of this spaces.

1)$\,$ For every $0< \|f\|_{p,\,\omega,\,\Omega}< \infty,$ $\displaystyle{I_{p,\,\omega}\left(\frac f{\|f\|_{p,\,\omega,\,\Omega}} \right)= 1}.$

If $\displaystyle{I_{p,\,\omega}\left(\frac f{\|f\|_{p,\,\omega,\,\Omega}} \right)< 1},$ we can find
$0< \lambda\le \|f\|_{p,\,\omega,\,\Omega}$ such that $\displaystyle{I_{p,\,\omega}\left(\frac f{\lambda} \right)< 1}.$ Indeed, let $\displaystyle{\lambda= \|f\|_{p,\,\omega,\,\Omega}\,I_{p,\,\omega}^{1/\overline p}\left(\frac f{\|f\|_{p,\,\omega,\,\Omega}}\right)}.$ Then $\lambda< \|f\|_{p,\,\omega,\,\Omega}$ and the inequality
$$
I_{p,\,\omega}\left(\frac f{\lambda} \right)= \int\limits_{\Omega} \left(\frac{|f(x)|\,\omega(x)}
{\|f\|_{p,\,\omega,\,\Omega}\,I_{p,\,\omega}^{1/\overline p}\left(\frac f{\|f\|_{p,\,\omega,\,\Omega}} \right)}\right)^{p(x)}\,dx
$$
$$
\le I_{p,\,\omega}^{-1}\left(\frac f{\|f\|_{p,\,\omega,\,\Omega}} \right)
\int\limits_{\Omega} \left(\frac{|f(x)|\,\omega(x)}{\|f\|_{p,\,\omega,\,\Omega}}\right)^{p(x)}\,dx
= 1
$$
is valid.
The obtained inequality contradicts to (2.1).

\begin{remark}
Note that  property 1) for non-weighted case was proved in [15].
\end{remark}

2)$\,$ $\displaystyle{\min\left\{\|f\|_{p,\,\omega,\,\Omega}^{\underline p},\,\|f\|_{p,\,\omega,\,\Omega}^{\overline p}\right)\le I_{p,\,\omega}(f)\le \max\left\{\|f\|_{p,\,\omega,\,\Omega}^{\underline p},\,\|f\|_{p,\,\omega,\,\Omega}^ {\overline p}\right)}.$

Let $\|f\|_{p,\,\omega,\,\Omega}\le 1.$  Using the property 1) we have
$$
I_{p,\,\omega}(f)=  \int\limits_{\Omega} \|f\|_{p,\,\omega,\,\Omega}^{p(x)}\left(\frac{|f(x)|\,\omega(x)} {\|f\|_{p,\,\omega,\,\Omega}}\right)^{p(x)}\,dx\le \|f\|_{p,\,\omega,\,\Omega}^{\underline p}\,\int\limits_{\Omega} \left(\frac{|f(x)|\,\omega(x)} {\|f\|_{p,\,\omega,\,\Omega}}\right)^{p(x)}\,dx= \|f\|_{p,\,\omega,\,\Omega}^ {\underline p}.
$$
Conversely, $I_{p,\,\omega}(f)\ge \|f\|_{p,\,\omega,\,\Omega}^ {\overline p}.$ Analogously, is consider the case
$\|f\|_{p,\,\omega,\,\Omega}\ge 1.$

3) The space $L_{p(x),\,\omega}(\Omega)$ is real linear spaces.

By using of the property 1), we have
$$
\int\limits_{\Omega} \left(\frac{|f(x)+ g(x)|\,\omega(x)}{2^{1/\underline p}\,\left(\|f\|_{p,\,\omega,\,\Omega}+ \|g\|_{p,\,\omega,\,\Omega}\right)}\right)^{p(x)}\,dx
$$
$$
\le \int\limits_{\Omega} \left(\frac{|f(x)|\,\omega(x)}{2^{1/\underline p}\,\left(\|f\|_{p,\,\omega,\,\Omega}+ \|g\|_{p,\,\omega,\,\Omega} \right)}\right)^{p(x)}\,dx+ \int\limits_{\Omega} \left(\frac{|g(x)|\,\omega(x)}{2^{1/\underline p}\,\left(\|f\|_{p,\,\omega,\,\Omega}+ \|g\|_{p,\,\omega,\,\Omega} \right)}\right)^{p(x)}\,dx
$$
$$
\le \int\limits_{\Omega} 2^{-\frac{p(x)}{\underline p}}\,\left(\frac{|f(x)|\,\omega(x)} {\|f\|_{p,\,\omega,\,\Omega}}\right)^{p(x)}\,dx+ \int\limits_{\Omega} 2^{-\frac{p(x)}{\underline p}}\,\left(\frac{g(x)|\,\omega(x)} {\|g\|_{p,\,\omega,\,\Omega}}\right)^{p(x)}\,\,dx
$$
$$
\le \frac 12\,\left(\int\limits_{\Omega} \left(\frac{|f(x)|\,\omega(x)} {\|f\|_{p,\,\omega,\,\Omega}}\right)^{p(x)}\,dx+ \int\limits_{\Omega} \left(\frac{|g(x)|\,\omega(x)} {\|g\|_{p,\,\omega,\,\Omega}}\right)^{p(x)}\,dx\right)= 1.
$$
Thus by Definition 1 $\displaystyle{\left\|f+ g\right\|_{p,\,\omega,\,\Omega}\le 2^{1/\underline p}\,\left(\|f\|_{p,\,\omega,\,\Omega}+ \|g\|_{p,\,\omega,\,\Omega}\right)}.$ Therefore $f+ g\in L_{p(x),\,\omega}(\Omega).$

Let $\alpha\in R\setminus \{0\}$ and $f\in L_{p(x),\,\omega}(\Omega).$
Now show that $\alpha f\in L_{p(x),\,\omega}(\Omega).$ We get
$$
\|\alpha\,f\|_{p,\,\omega,\,\Omega}= \inf\left\{\lambda> 0:\;\;\int\limits_{\Omega} \left(\frac{|\alpha\,f(x)|\,\omega(x)}{\lambda}\right)^{p(x)} \,dx\le 1\right\}
$$
$$
= \inf\left\{\lambda> 0:\;\;\int\limits_{\Omega} \left(\frac{|f(x)|\,\omega(x)}{\frac\lambda{|\alpha|}}
\right)^{p(x)} \,dx\le 1\right\}
$$
We substitute $\lambda= |\alpha|\,\mu.$ Then
$$
\inf\left\{\lambda> 0:\;\;\int\limits_{\Omega} \left(\frac{|f(x)|\,\omega(x)}{\frac\lambda{|\alpha|}}
\right)^{p(x)} \,dx\le 1\right\}
$$
$$
=\inf\left\{|\alpha|\mu> 0:\;\;\int\limits_{\Omega} \left(\frac{|f(x)|\,\omega(x)}{\mu}\right)^{p(x)} \,dx\le 1\right\}
$$
$$
=|\alpha|\, \inf\left\{\mu> 0:\;\;\int\limits_{\Omega} \left(\frac{|f(x)|\,\omega(x)}{\mu}\right)^{p(x)} \,dx\le 1\right\}= |\alpha|\,\|f\|_{p,\,\omega,\,\Omega}.
$$
For $f= 0$ this fact is trivially. Hence implies that the variable Lebesgue space $L_{p(x),\,\omega}(\Omega)$
is real linear space.

4) Let $\|f\|_{p,\,\omega,\,\Omega}= 0.$ Then we proved that $f= 0$ a.e. $x\in \Omega.$

If $\|f\|_{p,\,\omega,\,\Omega}= 0,$ then by (2.1) for all $\lambda> 0,$ $\displaystyle{I_{p,\,\omega}\left(\frac f{\lambda}\right)\le 1}.$ For any $\mu> 0$ and $\varepsilon\in (0,\,1),$ we have
$$
I_{p,\,\omega}\left(\frac f{\mu}\right)= \int\limits_{\Omega} \varepsilon^{p(x)}\left(\frac{|f(x)|\,\omega(x)}
{\varepsilon\,\mu}\right)^{p(x)} \,dx\le \varepsilon^{\underline p}\,I_{p,\,\omega}\left(\frac f{\varepsilon \mu}\right)\le \varepsilon^{\underline p}.
$$
Since $\varepsilon$ be any number from $(0, 1),$ then $\displaystyle{I_{p,\,\omega}\left(\frac f{\mu}\right)= 0}$ for all $\mu> 0.$ Therefore \linebreak $\displaystyle{\int\limits_{\Omega} \left(\frac{|f(x)|\,\omega(x)}
{\mu}\right)^{p(x)} \,dx=0}$ and thus $f= 0$ a.e. $x\in \Omega.$

5) Let $|f(x)|\le |g(x)|$ for a.e. $x\in \Omega.$ Then $\|f\|_{p,\,\omega,\,\Omega}\le \|g\|_{p,\,\omega,\,\Omega}.$

Indeed, by virtue of property 1) we have
$$
\int\limits_{\Omega} \left(\frac{|f(x)|\,\omega(x)}{\|g\|_{p,\,\omega,\,\Omega}}\right)^{p(x)}\,dx= \int\limits_{\Omega}\left(\frac{|f(x)|}{|g(x)|}\,\frac{|g(x)|\,\omega(x)}{\|g\|_{p,\,\omega,\,\Omega}}
\right)^{p(x)}\,dx
$$
$$
\le \int\limits_{\Omega}\left(\frac{|g(x)|\,\omega(x)}{\|g\|_{p,\,\omega,\,\Omega}}\right)^{p(x)}\,dx= 1.
$$
Thus by Definition 1 $\|f\|_{p,\,\omega,\,\Omega}\le \|g\|_{p,\,\omega,\,\Omega}.$

\begin{lemma}
Let $0< \underline p\le p(x)\le \overline p< 1$ and $f,g\in L_{p(x),\,\omega}(\Omega).$
Then
$$
\left\||f|+ |g|\right\|_{p,\,\omega,\,\Omega}\ge \|f\|_{p,\,\omega,\,\Omega}+ \|g\|_{p,\,\omega,\,\Omega}.
$$
\end{lemma}

{\bf Proof.} First we show that the function $h(t)= t^r,$ for $0< r< 1$ and $t> 0$ is concave. Let $\alpha+ \beta= 1,$ where $\alpha,\, \beta\ge 0.$ We proved that $\left(\alpha+ \beta\,t\right)^r\ge \alpha+ \beta\,t^r.$
We consider the function $\displaystyle{F(t)= \frac{\left(\alpha+ \beta\,t\right)^r}{\alpha+ \beta\,t^r}}.$
Differentiating by $t$ and after some calculation we have
$$
F'(t)= \frac{\alpha\,\beta \,p\left(\alpha+ \beta\,t\right)^{r- 1}
\left(1- t^{r- 1}\right)}{\left(\alpha+ \beta\,t^r\right)^2}.
$$
Since $r- 1< 0,$ then $t= 1$ is minimal value of the function $F$ for all $t> 0.$ Therefore $F(t)\ge F(1)= 1.$
Thus $\left(\alpha+ \beta\,t\right)^r\ge \alpha+ \beta\,t^r.$ Taking $\displaystyle{t= \frac{t_2}{t_1}}$ in last inequality we have $\displaystyle{\left(\alpha\,t_1+ \beta\,t_2\right)^r\ge \alpha\,t_1^r+ \beta\,t_2^r},$
i.e. the function $h(t)= t^r$ is concave.

Now we show a requiring inequality. It is obvious that the case $f=g= 0$ a.e. $x\in \Omega$ is trivial. Let $\|f\|_{p,\,\omega,\,\Omega}> 0$ and $\|g\|_{p,\,\omega,\,\Omega}> 0.$
Using concavity property of power function and property 1), we get
$$
I_{p,\,\omega}\left(\frac{|f|+ |g|}{\|f\|_{p,\,\omega,\,\Omega}+ \|g\|_{p,\,\omega,\,\Omega}}\right)= \int\limits_{\Omega} \left(\frac{|f(x)|+ |g(x)|}{\|f\|_{p,\,\omega,\,\Omega}+ \|g\|_{p,\,\omega,\,\Omega}}\,
\omega(x)\right)^{p(x)}\,dx=
$$
$$
\int\limits_{\Omega} \left(\frac{\|f\|_{p,\,\omega,\,\Omega} \frac{|f(x)|}{\|f\|_{p,\,\omega,\,\Omega}}+ \|g\|_{p,\,\omega,\,\Omega}\frac{|g(x)|} {\|g\|_{p,\,\omega,\,\Omega}}}{\|f\|_{p,\,\omega,\,\Omega}+ \|g\|_{p,\,\omega,\,\Omega}}\,\omega(x)\right)^{p(x)}\,dx
$$
$$
=\int\limits_{\Omega} \left(\frac{\|f\|_{p,\,\omega,\,\Omega}}{\|f\|_{p,\,\omega,\,\Omega}+ \|g\|_{p,\,\omega,\,\Omega}}\frac{|f(x)|}{\|f\|_{p,\,\omega,\,\Omega}}+ \frac{\|g\|_{p,\,\omega,\,\Omega}}{\|f\|_{p,\,\omega,\,\Omega}+ \|g\|_{p,\,\omega,\,\Omega}} \frac{|g(x)|}{\|g\|_{p,\,\omega,\,\Omega}}\right)^{p(x)}\,[\omega(x)]^{p(x)}\,dx
$$
$$
\ge \frac{\|f\|_{p,\omega,\Omega}}{\|f\|_{p,\omega,\Omega}+ \|g\|_{p,\omega,\Omega}}
\int\limits_{\Omega}\left(\frac{|f(x)|\omega(x)}{\|f\|_{p,\omega,\Omega}}\right)^{p(x)} dx+
\frac{\|g\|_{p,\omega,\Omega}}{\|f\|_{p,\omega,\Omega}+ \|g\|_{p,\omega,\Omega}}
\int\limits_{\Omega}\left(\frac{|g(x)|\omega(x)}{\|g\|_{p,\omega,\Omega}}\right)^{p(x)} dx
$$
$$
= \frac{\|f\|_{p,\omega,\Omega}}{\|f\|_{p,\omega,\Omega}+ \|g\|_{p,\omega,\Omega}}+
\frac{\|g\|_{p,\omega,\Omega}}{\|f\|_{p,\omega,\Omega}+ \|g\|_{p,\omega,\Omega}}= 1.
$$
Thus $\displaystyle{\left\||f|+ |g|\right\|_{p,\,\omega,\,\Omega}\ge \|f\|_{p,\,\omega,\,\Omega}+ \|g\|_{p,\,\omega,\,\Omega}}.$ In addition, note that the inequality in the form $\displaystyle{\left\|f+ g\right\|_{p,\,\omega,\,\Omega}\ge \|f\|_{p,\,\omega,\,\Omega}+ \|g\|_{p,\,\omega,\,\Omega}}$ doesn't holds for
$L_{p(x),\,\omega}(\Omega).$ Indeed, taking $g= -f$ we can see $0\ge 2\,\|f\|_{p,\,\omega,\,\Omega},$ which is valid
only for $f= 0$ a.e. $x\in \Omega.$

This proves the Lemma 1.

\begin{theorem} Let $0< \underline p\le p(x)\le \overline p< 1$ and $\displaystyle{p'(x)= \frac{p(x)}{p(x)- 1}}$
and $\omega$ be a weight function defined on $\Omega.$
Then the inequality
$$
\int\limits_{\Omega} |f(x)\,g(x)|\,dx\ge \left(\frac 1{\overline p}+ \frac 1{\overline p'}\right)
\|f\|_{p,\omega,\Omega}\,\|g\|_{p',\,\omega^{-1},\,\Omega} \eqno(2.2)
$$
holds for every $f\in L_{p(x),\omega}(\Omega),$ $g\in L_{p'(x),\omega^{-1}}(\Omega)$ and $0< |g(x)|< \infty.$
\end{theorem}

{\bf Proof.} We consider the function $\displaystyle{G(t)= \frac{t^s}{s}+ \frac{t^{-s'}}{s'}},$ where $t> 0,$
$0<s= const<1$ and $\displaystyle{s'= \frac{s}{s- 1}}.$ Differentiating by $t$ we have
$$
G'(t)= t^{s-1}- \frac 1{t^{s'+ 1}}= \frac{t^{ss'}- 1}{t^{s'+ 1}},
$$
where $s+ s'= s s'< 0.$ Therefore the point $t= 1$ is maximal value of the function $G(t)$ for all $t> 0.$
Thus $G(t)\le G(1)= 1,$ i.e., $\displaystyle{\frac{t^s}{s}+ \frac{t^{-s'}}{s'}\le 1}.$ If we take
$\displaystyle{t= \frac {a^{1/s'}}{b^{1/s}}},$ then
$$
ab\ge \frac{a^s}{s}+ \frac{b^{s'}}{s'}, \eqno(2.3)
$$
where $a, b> 0.$

Putting $\displaystyle{a= \frac{|f(x)|\omega(x)}{\|f\|_{p,\,\omega,\,\Omega}}},$ $\displaystyle{b= \frac{|g(x)|\omega^{-1}(x)}{\|g\|_{p',\,\omega^{-1},\,\Omega}}},$ $s= s(x)= p(x),$ $s'= s'(x)= p'(x)$
in inequality (2.3) and using the property 1) we have
$$
\int\limits_{\Omega} \frac{|f(x)\,g(x)|}{\|f\|_{p,\,\omega,\,\Omega} \|g\|_{p',\,\omega^{-1},\,\Omega}} dx\ge
\int\limits_{\Omega} \frac 1{p(x)} \left(\frac{|f(x)|\omega(x)}{\|f\|_{p,\,\omega,\,\Omega}}\right)^{p(x)}\,dx+
\int\limits_{\Omega} \frac 1{p'(x)} \left(\frac{|g(x)|\omega^{-1}(x)}{\|g\|_{p',\,\omega^{-1},\,\Omega}} \right)^{p'(x)}\,dx
$$
$$
\ge \frac 1{\overline p}\,\int\limits_{\Omega}\left(\frac{|f(x)|\omega(x)}{\|f\|_{p,\,\omega,\,\Omega}}\right) ^{p(x)}\,dx+ \frac 1{\overline p'}\,\int\limits_{\Omega} \left(\frac{|g(x)|\omega^{-1}(x)} {\|g\|_{p',\,\omega^{-1},\,\Omega}} \right)^{p'(x)}\,dx= \frac 1{\overline p}+ \frac 1{\overline p'}.
$$
Thus the inequality (2.2) is proved.

\begin{remark} Note that in the proof of Lemma 1, the expression $\|g\|_{p',\,\omega^{-1},\,\Omega}$ was used for negative values of the conjugate function. It should be understood as follows
$$
\|g\|_{p',\,\omega^{-1},\,\Omega}:= \inf\limits\left\{\lambda> 0: \;\;\int\limits_{\Omega} \left(\frac{\left[|g(x)|\,\omega^{-1}(x)\right]^{-1}}{\lambda^{-1}}\right)^{-p'(x)}\le 1\right\}
$$
$$
= \inf\limits\left\{\frac 1{\mu}> 0: \;\;\int\limits_{\Omega} \left(\frac{\left[|g(x)|\,\omega^{-1}(x)\right]^{-1}} {\mu}\right)^{-p'(x)}\le 1\right\}=
$$
$$
=\sup\limits\left\{\mu> 0: \;\;\int\limits_{\Omega} \left(\frac{|g(x)|\,\omega^{-1}(x)}{\mu}\right)^{p'(x)}\le 1\right\}.
$$
\end{remark}

\begin{theorem} Let $0< \underline p\le p(x)\le q(x)\le \overline q< 1$ and $\displaystyle{r(x)= \frac{p(x)\,q(x)}{q(x)- p(x)}}.$ Suppose that $\omega_1$ and $\omega_2$ are weights functions defined
in $\Omega$ and satisfying the condition
$$
\left\|\frac{\omega_1}{\omega_2} \right\|_{r,\,\Omega}< \infty.
$$
Then the inequality
$$
\|f\|_{p,\,\omega_1,\,\Omega}\le \left(A+ B+ \left\|\chi_{\Omega_2}\right\|_{L_{\infty} \left(\Omega\right)}\right)^{1/\underline p} \left\|\frac{\omega_1}{\omega_2} \right\|_{L_{r(\cdot)}\left(\Omega\right)}\,\|f\|_{q,\omega_2,\Omega},
$$
holds for every $f\in L_{q(x),\omega_2}(\Omega),$ where $\Omega_1= \left\{x\in \Omega:\,p(x)< q(x)\right\},$
$\Omega_2= \left\{x\in \Omega:\,p(x)= q(x)\right\}$ and $\displaystyle{A= \sup\limits_{x\in \Omega_1}
\frac{p(x)}{q(x)}},$  $\displaystyle{B= \sup\limits_{x\in \Omega_1} \frac{q(x)- p(x)}{q(x)}}\,$ and
$\,\displaystyle{\left\|\frac{\omega_1}{\omega_2}\right\|_{L_{r(\cdot)}\left(\Omega\right)}= \left\|\frac{\omega_1}{\omega_2} \right\|_{L_{r(\cdot)}\left(\Omega_1\right)}+ \left\|\frac{\omega_1}{\omega_2} \right\|_{L_{\infty} \left(\Omega_2\right)}}.$
\end{theorem}

{\bf Proof.} We have
$$
\|f\|_{p,\,\omega_1,\,\Omega_2}= \left\|f\,\omega_2\,\frac{\omega_1}{\omega_2}\right\|_{p,\,\omega_1,\,\Omega_2}\le
\left\|\frac{\omega_1}{\omega_2}\right\|_{L_{\infty}\left(\Omega_2\right)}\,\left\|f\omega_2\right\|_{p,\,\Omega_2}
$$
$$
=
\left\|\frac{\omega_1}{\omega_2}\right\|_{L_{\infty}\left(\Omega_2\right)}\,\left\|f\chi_{\Omega_2}\right\|_{p,\,
\omega_2,\,\Omega}\le\left\|\frac{\omega_1}{\omega_2}\right\|_{L_{\infty}\left(\Omega_2\right)}\,\left\|\chi_{\Omega_2}
\right\|_{L_{\infty} \left(\Omega\right)}\left\|f\right\|_{p,\,\omega_2,\,\Omega}.
$$
Therefore $\displaystyle{\left\|\frac{f}{\left\|\frac{\omega_1}{\omega_2}\right\|_{L_{\infty}\left(\Omega_2\right)}\,
\left\|f\right\|_{p,\,\omega_2,\,\Omega}}\right\|_{p,\,\omega_1,\,\Omega_2}\le \left\|\chi_{\Omega_2}
\right\|_{L_{\infty} \left(\Omega\right)}}\le 1.$ By virtue of property 1)
$$
\int\limits_{\Omega_2}\left(\frac{|f(x)|\,\omega_1(x)}{\left\|\frac{\omega_1}{\omega_2}\right\|_{L_{\infty}
\left(\Omega_2\right)}\,\|f\|_{p,\omega_2,\Omega}}\right)^{p(x)}\,dx\le \left\|\chi_{\Omega_2}
\right\|_{L_{\infty} \left(\Omega\right)}^{\underline p}= \left\|\chi_{\Omega_2}
\right\|_{L_{\infty} \left(\Omega\right)}. \eqno(2.4)
$$
It is well known that the inequality (2.3) for $s> 1$ is Young's inequality, i.e.
$$
ab\le \frac{a^s}{s}+ \frac{b^{s'}}{s'}, \eqno(2.5)
$$
where $\displaystyle{s'= \frac s{s- 1}}.$  We take $\displaystyle{s= s(x)= \frac{q(x)}{p(x)}},$ $\displaystyle{a= \left(\frac{|f(x)|\,\omega_2(x)}{\|f\|_{q,\omega_2,\Omega_1}}\right)^{p(x)}}$ and $\displaystyle{b= \left[\frac{\omega_1(x)} {\omega_2(x)}/\left\|\frac{\omega_1}{\omega_2}\right\|_{r,\,\Omega_1}\right]^{p(x)}}.$
Thus $\displaystyle{s'= s'(x)= \frac{q(x)}{q(x)- p(x)}}$ and from inequality (2.5), we have
$$
\left(\frac{|f(x)|\,\omega_1(x)}{\left\|\frac{\omega_1}{\omega_2}\right\|_{r,\,\Omega_1}\,
\|f\|_{q,\omega_2,\Omega_1}}\right)^{p(x)}\le \frac{p(x)}{q(x)}\,\left(\frac{|f(x)|\,\omega_2(x)} {\|f\|_{q,\omega_2,\Omega_1}}\right)^{q(x)}+ \frac{q(x)- p(x)}{q(x)}\,\left[\frac{\frac{\omega_1(x)} {\omega_2(x)}}{\left\|\frac{\omega_1}{\omega_2}\right\|_{r,\,\Omega_1}}\right]^{r(x)}
$$
$$
\le A\,\left(\frac{|f(x)|\,\omega_2(x)} {\|f\|_{q,\omega_2,\Omega_1}}\right)^{q(x)}+
B\,\left[\frac{\frac{\omega_1(x)} {\omega_2(x)}}{\left\|\frac{\omega_1}{\omega_2}
\right\|_{r,\,\Omega_1}}\right]^{r(x)}.
$$
Obviously, $1\le A+ B\le 2.$
Integrating by $\Omega_1,$ using the property 1), we get
$$
\int\limits_{\Omega_1}\left(\frac{|f(x)|\,\omega_1(x)}{\left\|\frac{\omega_1}{\omega_2}\right\|_{r,\,\Omega_1}\,
\|f\|_{q,\omega_2,\Omega_1}}\right)^{p(x)}\,dx
$$
$$
\le A\,\int\limits_{\Omega_1}\left(\frac{|f(x)|\,\omega_2(x)} {\|f\|_{q,\omega_2,\Omega_1}}\right)^{q(x)}\,dx +
B\,\int\limits_{\Omega_1}\left[\frac{\frac{\omega_1(x)} {\omega_2(x)}}{\left\|\frac{\omega_1}{\omega_2}
\right\|_{r,\,\Omega_1}}\right]^{r(x)}\,dx \le A+ B. \eqno(2.6)
$$
From (2.4) and (2.6) implies that
$$
\int\limits_{\Omega}\left(\frac{|f(x)|\,\omega_1(x)}{\left\|\frac{\omega_1}{\omega_2}\right\|_
{L_{r(\cdot)}\left(\Omega\right)}\,\|f\|_{q,\omega_2,\Omega}}\right)^{p(x)}\,dx= \int\limits_{\Omega_1}\left(\frac{|f(x)|\,\omega_1(x)}{\left\|\frac{\omega_1}{\omega_2}\right\|_
{L_{r(\cdot)}\left(\Omega\right)}\,\|f\|_{q,\omega_2,\Omega}}\right)^{p(x)}\,dx
$$
$$
+\int\limits_{\Omega_2}\left(\frac{|f(x)|\,\omega_1(x)}{\left\|\frac{\omega_1}{\omega_2}\right\|_
{L_{r(\cdot)}\left(\Omega\right)}\,\|f\|_{q,\omega_2,\Omega}}\right)^{p(x)}\,dx \le \int\limits_{\Omega_1}\left(\frac{|f(x)|\,\omega_1(x)}{\left\|\frac{\omega_1}{\omega_2}\right\|_
{L_{r(\cdot)}\left(\Omega_1\right)}\,\|f\|_{q,\omega_2,\Omega_1}}\right)^{p(x)}\,dx
$$
$$
+ \int\limits_{\Omega_2}\left(\frac{|f(x)|\,\omega_1(x)}{\left\|\frac{\omega_1}{\omega_2}\right\|_
{L_{\infty}\left(\Omega_2\right)}\,\|f\|_{q,\omega_2,\Omega}}\right)^{p(x)}\,dx\le A+ B+ \left\|\chi_{\Omega_2}\right\|_{L_{\infty}\left(\Omega\right)}.
$$
From last inequality we have
$$
1\ge \int\limits_{\Omega}\left(\frac{|f(x)|\,\omega_1(x)}{\left(A+ B+ \left\|\chi_{\Omega_2}\right\|_{L_{\infty} \left(\Omega\right)}\right)^{1/p(x)} \left\|\frac{\omega_1}{\omega_2} \right\|_{L_{r(\cdot)} \Omega} \,\|f\|_{q,\omega_2,\Omega}}\right)^{p(x)}\,dx
$$
$$
\ge \int\limits_{\Omega}\left(\frac{|f(x)|\,\omega_1(x)}{\left(A+ B+ \left\|\chi_{\Omega_2}\right\|_{L_{\infty} \left(\Omega\right)}\right)^{1/\underline p} \left\|\frac{\omega_1}{\omega_2} \right\|_{L_{r(\cdot)} \Omega} \,\|f\|_{q,\omega_2,\Omega}}\right)^{p(x)}\,dx.
$$
Thus
$$
\left\|f\right\|_{p,\,\omega_1,\,\Omega}\le \left(A+ B+ \left\|\chi_{\Omega_2}\right\|_{L_{\infty} \left(\Omega\right)}\right)^{1/\underline p} \left\|\frac{\omega_1}{\omega_2} \right\|_{L_{r(\cdot)}\left(\Omega\right)}\,\|f\|_{q,\omega_2,\Omega}.
$$

The theorem is proved.
\begin{remark} Note that Theorem 2 in the case $\omega_1= \omega_2= 1$ and $|\Omega|< \infty$ was proved in [15].
In the case $1\le \underline p\le p(x)\le q(x)\le \overline q< \infty$ for general measures Theorem 2 was proved in [4].
\end{remark}

The following theorems are known.
\begin{theorem}{[1]} Let $1\le \underline p\le p(x)\le q(y)\le \overline q<
\infty$ for all $x\in \Omega_1\subset R^n$ and $y\in \Omega_2\subset
R^m.$  If $p(x)\in C\left(\Omega_1\right),$ then the inequality
$$
\left\|\|f\|_{L_{p(\cdot)}\left(\Omega_1\right)}\right\|_{L_{q(\cdot)}\left(\Omega_2\right)}\le
C_{p,q}\,\left\|\|f\|_{L_{q(\cdot)}\left(\Omega_2\right)}\right\|_
{L_{p(\cdot)}\left(\Omega_1\right)}
$$
is valid, where $\displaystyle{C_{p, q}= \left(\left\|\chi_{\Delta_1}\right\|_\infty+ \left\|\chi_{\Delta_2}\right\|_\infty+ \frac{\overline p}{\underline q}- \frac{\underline p}{\overline q}\right)\left(\left\|\chi_{\Delta_1}\right\|_\infty+ \left\|\chi_{\Delta_2}\right\|_\infty\right)},$
$\underline q= \mbox {ess}\,\inf\limits_{\Omega_2} q(x),$ $\overline q= \mbox {ess}\,\sup\limits_{\Omega_2} q(x),$ $\Delta_1=\left\{(x, y)\in\Omega_1\times \Omega_2:\, p(x)= q(y)\right\},$ $\Delta_2= \Omega_1\times \Omega_2\setminus \Delta_1$ and
$C\left(\Omega_1\right)$ is the space of continuous functions in
$\Omega_1$ and $f:\Omega_1\times \Omega_2\rightarrow R$ is any
measurable function such that
$$
\left\|\|f\|_{q,\Omega_2}\right\|_{p,\Omega_1}= \inf \left\{ \mu>
0:\;\;\int\limits_{\Omega_1}\left(\frac{\|f(x,\cdot)\|_{q(\cdot),\Omega_2}}
{\mu}\right)^{p(x)}\,dx\le 1\right\}< \infty.
$$
\end{theorem}

The following lemmas are known.
\begin{lemma}{[6]} Let $0< s< 1,$  $-\infty< a< b\le \infty$ and $f$ is non-negative and
decreasing function defined on $(a, b).$ Then
$$
\left(\int\limits_{a}^{b} f(x)\,dx\right)^{s}\le s\,\int\limits_{a}^{b} f^s(x)\,(x- a)^{s- 1}\,dx.
$$
\end{lemma}

\begin{lemma}{[6]} Let $0< s< 1,$  $-\infty\le a< b< \infty$ and $f$ is non-negative and
increasing function defined on $(a, b).$ Then
$$
\left(\int\limits_{a}^{b} f(x)\,dx\right)^{s}\le s\,\int\limits_{a}^{b} f^s(x)\,(b- x)^{s- 1}\,dx.
$$
\end{lemma}

\vspace{3mm}
\begin{center}
3. {\bf On a topology of the spaces $L_{p(x),\,\omega}$ for $0< p(x)< 1$}
\end{center}

Now we formulate some definitions which be characterized of the topology in general vector spaces.

\begin{definition} A subset $G$ of a vector space X is called convex if, for any $x_1, x_2,\ldots,x_m\in G,$
$\displaystyle{\sum\limits_{i= 1}^m \alpha_i\,x_i\in G},$ where $\displaystyle{\sum\limits_{i= 1}^m \alpha_i= 1}$
and $\alpha_i\ge 0,$ $i= 1,2,\ldots,m.$ In particular, the subset contains the average $\frac 1m\,\sum\limits_
{i= 1}^m x_i.$
\end{definition}

\begin{definition} A topological vector space $X$ is called locally convex if the convex open sets are a base
for the topology, i.e., any open set $U\subset X$ around a point, there is a convex open set $C$ containing that
point such $C\subset X.$
\end{definition}

We show that the weighted variable Lebesgue spaces $L_{p(x),\,\omega}(\Omega)$ isn't locally convex.

\begin{lemma} Let $0< \underline p\le p(x)\le q(x)\le \overline q< 1$ and $\omega$ be a weight function defined on
$\Omega$ and $0< \omega(x)< \infty$ a.e. $x\in \Omega.$ Then weighted variable Lebesgue spaces $L_{p(x),\,\omega} (\Omega)$ isn't locally convex.
\end{lemma}

{\bf Proof.} It is obvious that $\displaystyle{\rho(f, g)= \int\limits_{\Omega} \left[\left|f(x)- g(x)\right|\omega(x)\right]^{p(x)}\,dx}$ is defined a metric on $L_{p(x),\omega}(\Omega).$
We consider any open ball neighborhoods $0:$
$$
U_{R}(0)= \left\{f\in L_{p(x),\,\omega}(\Omega):\; \rho(f,0)= I_{p(x),\omega}(f)< R\right\}.
$$
We will show that, for any $\varepsilon> 0,$ the $\varepsilon-$ball neighborhoods zero contains functions whose
average lies outside the ball of radius $R.$

Suppose $\varepsilon> 0$ and $m\ge 1.$ We select $m$ disjoint intervals $A_1, A_2,\ldots,A_m$ in $\Omega,$ which
need not cover of all $\Omega.$ We put $\displaystyle{f_k= \left(\frac{\varepsilon}
{\omega\left(A_k\right)} \right)^{1/p(x)}\chi_{A_k}},$ where $\displaystyle{\omega\left(A_k\right)= \int\limits_{A_k}\left[\omega(x)\right]^{p(x)} dx}$ and $k= 1,2,\ldots,m.$ Then $\displaystyle{I_{p,\omega} \left(f_k\right)= \frac{\varepsilon}{\omega\left(A_k\right)}\int\limits_{A_k} \left[\omega(x)\right]^
{p(x)} dx= \varepsilon},$ and so every $f_k$ is at distance $\varepsilon$ from $0.$ But, since the functions
$f_k$ are supported on disjoint sets, their average $\displaystyle{g_m= \frac 1m\,\sum\limits_{i= 1}^m f_i}$ satisfies
$$
I_{p,\omega} \left(g_m\right)= \int\limits_{\Omega} g_m^{p(x)}(x)\,dx=
\int\limits_{\Omega} \frac 1{m^{p(x)}} \left(\sum\limits_{i= 1}^m f_i\right)^{p(x)}
\left[\omega(x)\right]^{p(x)}\,dx
$$
$$
\ge \frac 1{m^{\overline p}}\,\sum\limits_{i= 1}^m\int\limits_{\Omega} \left(f_i(x)\,\omega(x)\right)^{p(x)}\,
dx= \frac{\varepsilon}{m^{\overline p}}\,\sum\limits_{i= 1}^m 1= m^{1- \overline p}\,\varepsilon.
$$
Then $I_{p,\omega} \left(g_m\right)\to \infty,$ for $m\to \infty$ (depending on $\varepsilon$).
Therefore $\rho(g_m, 0)\to \infty,$ for $m\to \infty.$ Thus the distance between $g_n$ and $0$ can be made
as large as desired.

The Lemma 4 is proved.

\begin{theorem} Let $0< \underline p\le p(x)\le q(x)\le \overline q< 1$ and $\omega$ be a weight function defined
on $\Omega$ and $0< \omega(x)< \infty$ a.e. $x\in \Omega.$ Then $\left[L_{p(x),\omega}(\Omega)\right]^
{\star}= \{0\},$ where $\star$ - be denoted dual space of $L_{p(x),\omega}(\Omega),$ i.e., is the space of continuous linear functionals from $L_{p(x),\omega}(\Omega)$ to $R.$
\end{theorem}

{\bf Proof.} We argue by contradiction. Let $\varphi\ne 0$ and $\varphi\in \left[L_{p(x),\omega}(\Omega)
\right]^{\star}.$ Let $\widetilde B(0, t)= \Omega\bigcap B(0,t),$ where $0< t< \infty.$

Suppose that $\varphi$ is linear continuous functional defined in $L_{p(x),\omega}(\Omega).$
Then we can find an $f\in L_{p(x),\,\omega}(\Omega)$ such that $\varphi(f)=1$.
Now, the map $t\mapsto f\chi_{\widetilde B(0, t)}$ is continuous, since
$|f|^{p(x)}\,\omega(x)$ is integrable:
$$
\int\limits_{\widetilde B(0, t_2)} |f(x)|^{p(x)}\,\omega(x)\,dx- \int\limits_{\widetilde B(0, t_1)} |f(x)|^{p(x)}\,\omega(x)\,dx= \int\limits_{\Omega\bigcap B_{t_1, t_2}}|f(x)|^p\,\omega(x)\,dx,\;\;\mbox{for }\;\;
t_1< t_2,
$$
where $B_{t_1, t_2}= \left\{x:\;\;t_1\le |y|< t_2\right\}.$
Thus we may choose $t\in (t_1, \infty)$ such that $\varphi\left(f\chi_{\widetilde B(0, t)}\right)= \varphi\left(f\chi_{\Omega\setminus \widetilde B(0, t)}\right)= \frac 12.$ Next, notice that $g=
f\chi_{\widetilde B(0, t)}$ and $h= f\chi_{\Omega\setminus \widetilde B(0, t)}$ satisfy
$$
\int\limits_{\Omega} |f(x)|^{p(x)}\,\omega(x)\,dx
= \int\limits_{\widetilde B(0, t)} |f(x)|^{p(x)}\,\omega(x)\,dx+
\int\limits_{\Omega\setminus\widetilde B(0, t)} |f(x)|^{p(x)}\,\omega(x)\,dx
= I_{p,\,\omega} \left(g\right)+ I_{p,\,\omega} \left(h\right).
$$
Thus, at least one of $I_{p,\,\omega} \left(g\right)$ or
$I_{p,\,\omega} \left(h\right)$ is less than $\displaystyle{\frac 12\,I_{p,\,\omega} \left(f\right)}.$
Let's say that $\displaystyle{I_{p,\,\omega} \left(g\right)\le \frac 12\,I_{p,\,\omega} \left(f\right)}.$
Then, $f_1= 2g$ satisfies
$$
\varphi(f_1)=1\quad\hbox{and}\quad
I_{p,\,\omega} \left(f_1\right)\le 2^{\overline p}\,I_{p,\,\omega} \left(g\right)\le
2^{\overline p- 1}\,I_{p,\,\omega} \left(f\right).
$$
By induction, we can find a sequence $\left\{f_n\right\}_{n\ge 1}$ in $L_{p(x),\,\omega}(\Omega)$ with
$$
\varphi\left(f_n\right)=1\quad\hbox{and}\quad
I_{p,\,\omega} \left(f_n\right)\le 2^{n(\overline p- 1)}\,I_{p,\,\omega} \left(f\right).
$$
It is obvious that $\overline p-1<0$ and $f_n\to 0$ in $L_{p(x),\,\omega}(\Omega)$ while $T(f_n)=1$.
Thus, $T=0$ is the only continuous linear functional.

\begin{theorem} Let $0< \underline p\le p(x)\le q(x)\le \overline q< 1$ and $\omega$ be
a weight function defined on $\Omega$ and $0< \omega(x)< \infty$ a.e. $x\in \Omega.$
Then the spaces $L_{p(x), \,\omega}(\Omega)$ is complete.
\end{theorem}

{\bf Proof.} Let $\left\{f_n\right\},$ $n\in N$ be a sequence in $L_{p(x),\,\omega}(\Omega)$
such that
$$
\left\|f_n- f_m\right\|_{p,\,\omega,\,\Omega}\to 0, \quad \mbox{for}\quad n, m\to \infty.
$$
From properties 1) implies that
$$
\int\limits_{\Omega}\left(\left|f_{n}- f_{m}\right|\,\omega(x)\right)^{p(x)}\,dx\to 0, \quad \mbox{for}\quad n, m\to \infty.
$$
We choose the subsequence $\{n_k\}$ such that
$$
A= \sum\limits_{k= 1}^{\infty} \int\limits_{\Omega} \left(\left|f_{n_{k+ 1}}- f_{n_k}\right|\,\omega(x)\right)^{p(x)}\,dx< \infty.
$$
Then for any $\ell\in N$
$$
\int\limits_{\Omega}\left[\sum\limits_{k= 1}^{\ell} \left(\left|f_{n_{k+ 1}}- f_{n_k}\right|\,\omega(x)
\right)\right]^{p(x)}\,dx\le \sum\limits_{k= 1}^{\ell} \int\limits_{\Omega} \left(\left|f_{n_{k+ 1}}- f_{n_k}\right|\,\omega(x)\right)^{p(x)}\,dx\le A.
$$
If $\ell\to \infty,$ then by monotone convergence theorem
$$
\int\limits_{\Omega}\left[\sum\limits_{k= 1}^{\infty} \left(\left|f_{n_{k+ 1}}- f_{n_k}\right|\,\omega(x) \right)\right]^{p(x)}\,dx\le A.
$$
Therefore,
$$
\sum\limits_{k= 1}^{\infty}\left|f_{n_{k+ 1}}- f_{n_k}\right|\,\omega(x)< \infty, \quad a.e.\quad x\in \Omega.
$$
Hence, by completeness of R, $f_{n_k}$ converges a.e. $x\in \Omega.$ We define a measurable function f
by
$$
f(x)= \left\{
            \begin{array}{l}
            \lim\limits_{k\to \infty} f_{n_k},\qquad\qquad\, for\;\; a.e.\;\; x\in \Omega\\
               0,\qquad\qquad\qquad\quad otherwise.\\
            \end{array}
            \right.
$$
Since $\displaystyle{\int\limits_{\Omega}\left(\left|f_{n}- f_{m}\right|\,\omega(x)\right)^{p(x)}\,dx\to 0, \quad \mbox{for}\quad n, m\to \infty},$ then $\displaystyle{\left|f_{n}- f_{m}\right|^{p(x)}\to 0},$ $n, m\to \infty.$
Given $\varepsilon> 0$ we can find $N_\varepsilon$ so that $n\ge N_\varepsilon$ implies
$$
\left|\left\{x:\;\;|f_n(x)- f_m(x)|^{p(x)}\right\}\right|= \int\limits_{\left\{x:\;\;|f_n(x)- f_m(x)|^{p(x)}\right\}} dx\le  \varepsilon,\quad \mbox{for}\quad m\ge n.
$$
In particular, $\displaystyle{\left|\left\{x:\;\;|f_n(x)- f_{n_k}(x)|^{p(x)}\right\}\right|\le  \varepsilon,\quad \mbox{for}\quad k\to \infty}.$ Hence, by Fatou's lemma for $n\ge N_\varepsilon,$ we have
$$
\left|\left\{x:\;\;|f_n(x)- f(x)|^{p(x)}\right\}\right|= \left|\lim\limits_{k\to \infty} inf\left\{x:\;\;|f_n(x)- f_{n_k}(x)|^{p(x)}\right\}\right|
$$
$$
\le\lim\limits_{k\to \infty} inf \left|\left\{x:\;\;|f_n(x)- f_{n_k}(x)|^{p(x)}\right\}\right|\le \varepsilon.
$$
Hence $f\in L_{p(x),\, \omega}(\Omega)$ and $\displaystyle{\int\limits_{\Omega}\left|\left(f_{n}- f\right)\,\omega(x)\right|^{p(x)}\,dx\to 0},$ for $n\to \infty.$

This completes the proof of Theorem 5.

\begin{remark}
Note that from property 5) and  Theorem 5 implies that the spaces $L_{p(x),\, \omega}(\Omega)$ is ideal.
\end{remark}

\vspace{5mm}
\begin{center}
4. {\bf Main results.}
\end{center}

We consider the classical Hardy operator and it's dual operator defined as
$$
Hf(x)= \frac 1x \int\limits_{0}^x f(t)\,dt,\;\; H^*f(x)= \frac 1x \int\limits_x^{\infty} f(t)\,dt
$$
where $f$ is nonnegative function on $(0, \infty).$

\begin{lemma} Let $0< \underline p\le p_n\le  \overline p \le 1,$ $p_n\ge p_{n+ 1}$ and $\left\{x_n\right\}_{n\ge 1}$ be any non-negative sequence of real numbers such that $x_n^{p_n}\ge x_{n+1}^{p_{n+ 1}}$ for any $n\in \Bbb{N}.$

Then
$$
\left(\sum\limits_{n= 1}^\infty x_n^{\frac{p_n}{\underline p}} \right)^{\underline p}\le
\sum\limits_{n= 1}^\infty x_n^{p_n} \left[n^{p_n}- (n- 1)^{p_n}\right]\le \sum\limits_{n= 1}^\infty x_n^{p_n}.\eqno(4.1)
$$
\end{lemma}

{\bf Proof.} First we proved that
$$
\left(\sum\limits_{n= 1}^m x_n^{\frac{p_n}{p_m}}\right)^{p_m}\le\sum\limits_{n= 1}^m x_n^{p_n}
\left[n^{p_n}- (n- 1)^{p_n}\right]. \eqno(4.2)
$$
We consider the function $\displaystyle{h(t)= \frac{(1+t)^q- 1}{t^q}},$ where $t\ge 0$ and $0< q< 1.$ It is obvious that $\displaystyle{h'(t)= \frac{q\,\left[1- (1+ t)^{q-1}\right]}{t^{q+ 1}}\ge 0}$ for all $t\ge 0.$ In particular, the function $h(t)$ is monotone increasing in the segment $[0, B].$ Therefore $h(t)\le h(B),$ i.e.,
$$
(1+t)^q\le 1+ t^q\left[\left(B^{-1}+ 1\right)^q- B^{-q}\right] \;\;\mbox{for\; any}\;\;\;0\le t\le B. \eqno(4.3)
$$
Since $x_1^{p_1}\ge x_2^{p_2},$ then $\displaystyle{x_2\le x_1^{\frac{p_1}{p_2}}}.$ Therefore taking $\displaystyle{t= \frac{x_2}{x_1^{\frac{p_1}{p_2}}}},$ $B= 1$ and $q= p_2$ in (4.3), we have
$$
\left(x_1^{\frac{p_1}{p_2}}+ x_2\right)^{p_2}\le x_1^{p_1}+ x_2^{p_2}\left(2^{p_2}- 1\right). \eqno(4.4)
$$
It is obvious that the inequality (4.4) be inequality (4.2) for $m= 2.$ By the condition of Lemma 2 $\,p_2\ge p_3$ and so $2^{p_3}\le 2^{p_2}.$ Since $x_3\le \frac{x_1^{\frac{p_1}{p_3}}+ x_2^{\frac{p_2}{p_3}}}2$ from (4.3) and (4.4) for $\displaystyle{t= \frac{x_3}{x_1^{\frac{p_1}{p_3}}+ x_2^{\frac{p_2}{p_3}}}},$ $\displaystyle{B= \frac 12}$ and $q= p_3,$ we get
$$
\left(x_1^{\frac{p_1}{p_3}}+ x_2^{\frac{p_2}{p_3}}+ x_3\right)^{p_3}\le \left(x_1^{\frac{p_1}{p_3}}+ x_2^{\frac{p_2}{p_3}}\right)^{p_3}+ x_3^{p_3}\left(3^{p_3}- 2^{p_3}\right)
$$
$$
\le x_1^{p_1}+ x_2^{p_2}\left(2^{p_3}- 1\right)+ x_3^{p_3}\left(3^{p_3}- 2^{p_3}\right)\le
x_1^{p_1}+ x_2^{p_2}\left(2^{p_2}- 1\right)+ x_3^{p_3}\left(3^{p_3}- 2^{p_3}\right).
$$
The last inequality is (4.1) for $m= 3.$ Clearly $x_1^{\frac{p_1}{p_{m+ 1}}}+ x_2^{\frac{p_2}{p_{m+ 1}}}+\ldots+ x_m^{\frac{p_m}{p_{m+ 1}}}+ x_{m+ 1}\ge (m+ 1)x_{m+ 1}.$ Hence $\displaystyle{x_{m+ 1}\le \frac{x_1^{\frac{p_1}{p_{m+ 1}}}+ x_2^{\frac{p_2}{p_{m+ 1}}}+\ldots+ x_m^{\frac{p_m}{p_{m+ 1}}}}{m}}.$ Therefore taking
$$
t= \frac{x_{m+ 1}}{x_1^{\frac{p_1}{p_{m+ 1}}}+ x_2^{\frac{p_2}{p_{m+ 1}}}+\ldots+ x_m^{\frac{p_m}{p_{m+ 1}}}},\, B= \frac 1m\;\,  \mbox{and}\;\; q= p_{m+ 1}
$$
in (4.3), we have
$$
\left(\sum\limits_{n= 1}^{m+ 1} x_n^{\frac{p_n}{p_{m+ 1}}}\right)^{p_{m+ 1}}= \left(\sum\limits_{n= 1}^m x_n^{\frac{p_n}{p_{m+ 1}}}+ x_{m+ 1}\right)^{p_{m+ 1}}\le
$$
$$
\left(\sum\limits_{n= 1}^m x_n^{\frac{p_n}{p_{m+ 1}}} \right)^{p_{m+ 1}}+ x_{m+ 1}^{p_{m+ 1}}\left[(m+ 1)^{p_{m+ 1}}- m^{p_{m+ 1}}\right]\le
$$
$$
\sum\limits_{n= 1}^m x_n^{p_n}\left[n^{p_n}- (n- 1)^{p_n}\right]+ x_{m+ 1}^{p_{m+ 1}}\left[(m+ 1)^{p_{m+ 1}}- m^{p_{m+ 1}}\right]=
$$
$$
\sum\limits_{n= 1}^{m+ 1} x_n^{p_n}\left[n^{p_n}- (n- 1)^{p_n}\right].
$$
By the induction principle the inequality (4.2) is proved for any $m\in \Bbb N.$

Since the sequence $\left\{p_n\right\}_{n\ge 1}$ is decreasing, then $\lim\limits_{n\to \infty} p_n= \underline p.$
Therefore passing to the limit at $m\to \infty$ in (4.2) we have the left part of inequality (4.1).
By using the inequality $n^{p_n}\le (n- 1)^{p_n}+ 1,$ we have the right part of inequality (4.1).

The Lemma 2 is proved.

{\bf Example 4.1.}  Let $\displaystyle{x_n=\left\{
            \begin{array}{l}
             n^{-\frac {\underline p}{2\,p_n}},\quad \mbox{for}\;\, n= k^2\\
            0,\qquad\quad\, \mbox{for}\;\, n\ne k^2,\\
            \end{array}
            \right.}$ and $\displaystyle{\overline p< \frac{\underline p+ 1}{2}}.$

It is obvious that the sequence $\left\{x_n^{p_n}\right\}_{n\ge 1}$ isn't monotone and $\displaystyle{\sum\limits_{n= 1}^\infty x_n^{\frac{p_n}{\underline p}}= \sum\limits_{k= 1}^\infty \frac 1{k}= +\infty}.$ On the other hand  $\;\displaystyle{n^{p_n}- (n- 1)^{p_n}\sim p_n\,n^{p_n- 1}\sim n^{p_n- 1}}$ for $n\to \infty.$ Therefore
$$
\sum\limits_{n= 1}^{\infty} x_n^{p_n} \left[n^{p_n}- (n- 1)^{p_n}\right]\sim
\sum\limits_{n= 1}^{\infty} x_n^{p_n}\,n^{p_n- 1}= \sum\limits_{k= 1}^{\infty}
k^{-\underline p+ 2\,p_k- 2}\le \sum\limits_{k= 1}^{\infty}
k^{2\,\overline p-\underline p - 2}.
$$
It is well known that the series $\displaystyle{\sum\limits_{k= 1}^{\infty}
k^{2\,\overline p-\underline p - 2}}$ is converges if and only if $\displaystyle{\overline p< \frac{\underline p+ 1}{2}}.$  Thus for $\displaystyle{\overline p< \frac{\underline p+ 1}{2}}$ the inequality (3.1) isn't holds.

The example show that the condition of monotonicity of sequence $\left\{x_n^{p_n}\right\}_{n\ge 1}$ is essential.

\begin{remark}
Note that Lemma 5 in the case $p_1= p_2=\ldots= p_n=\ldots= p= const$ was proved in [5].
\end{remark}

\begin{theorem} Let $x\in (0, \infty),$ $0< \underline p\le p(x)\le q(x)\le  \overline q< 1,$ $\displaystyle{r(x)= \frac{\underline p\,p(x)}{p(x)- \underline p}}$ and $f(x)$ are non-negative and decreasing function defined on $(0, \infty).$ Suppose $\omega_1$ and $\omega_2$ are weight functions defined on $(0, \infty).$

Then for any $f\in L_{p(x),\,\omega_1}(0,\,\infty)$ the inequality
$$
\left\|Hf\right\|_{L_{q(\cdot),\,\omega_2}(0, \infty)}\le {\underline p}^{\frac 1{\underline p}}\,c_{p,q}\,d_p\, \left\|\frac{t^{1/p'}\,\left\|\frac{\omega_2}{x}\right\|_{L_{q(\cdot)}(t, \infty)}}{\omega_1} \right\|_{L_{r(\cdot)}(0, \infty)}\,\left\|f\right\|_{L_{p(\cdot),\,\omega_1}(0, \infty)},
$$
where $\displaystyle{c_{p,q}= \left(\left\|\chi_{\Delta_1}\right\|_{L_{\infty}(0,\,\infty)}+ \left\|\chi_{\Delta_2}\right\|_{L_{\infty}(0, \infty)}+ \underline p\left(\frac 1{\underline q}- \frac 1{\overline q}\right)\right)\left(\left\|\chi_{S_1}\right\|_{L_\infty (0, \infty)}+ \left\|\chi_{S_2}\right\|_{L_\infty (0, \infty)} \right)},$ $S_1= \left\{x\in (0, \infty):\; p(x)= \underline p\right\},$ $S_2= (0, \infty)\setminus S_1,$
and $\,\displaystyle{d_p=\left(1+ \frac{\overline p- \underline p}{\overline p}+ \left\|\chi_{S_1}\right\|_{L_{\infty}(0,\,\infty)}\right)^{1/\underline p}}.$
\end{theorem}

{\bf Proof.} Taking $a= 0,$ $b= x$ and $s= \underline p$ and apply Lemma 2 and property 5), we have
$$
\left\|Hf\right\|_{L_{q(\cdot),\,\omega_2}(0, \infty)}= \left\|\omega_2 Hf\right\|_{L_{q(\cdot)}(0, \infty)}= \left\|\frac{\omega_2}{x} \int\limits_{0}^x f(t)\,dt\right\|_{L_{q(\cdot)}(0, \infty)}
$$
$$
\le {\underline p}^{\frac 1{\underline p}}\left\|\frac{\omega_2(x)}{x} \left(\int\limits_{0}^x f^{\underline p}(t)\,t^{\underline p- 1}\,dt\right)^{1/\underline p}\right\|_{L_{q(\cdot)}(0, \infty)}.
$$
Now applied Theorem 3, we get
$$
\left\|\frac{\omega_2(x)}{x}\,\left(\int\limits_{0}^x f^{\underline p}(t)\,t^{\underline p- 1}\,dt\right)^{1/\underline p}\right\|_{L_{q(\cdot)}(0, \infty)}
$$
$$
= \left\|\left(\int\limits_{0}^\infty f^{\underline p}(t)\,\chi_{(0, \,x)}(t)\,\left[\frac{\omega_2(x)}{x}\right]^{\underline p}\,t^{\underline p- 1}\,dt\right)^{1/\underline p}\right\|_{L_{q(\cdot)}(0, \infty)}
$$
$$
= \left\|\int\limits_{0}^\infty f^{\underline p}(t)\,\chi_{(0,\,x)}(t)\,\left[\frac{\omega_2(x)}{x}\right]^
{\underline p}\,t^{\underline p- 1}\,dt\right\|_{L_{\frac{q(\cdot)}{\underline p}}(0, \infty)}^{1/\underline p}
$$
$$
\le c_{p,q}\,\left(\int\limits_{0}^\infty \left\| f^{\underline p}(t)\,\chi_{(0, \,x)}(t)\,\left[\frac{\omega_2(x)}{x}\right]^{\underline p}\,
\,t^{\underline p- 1}\right\|_{L_{\frac{q(\cdot)}{\underline p}}(0, \infty)}\,dt\right)^{1/\underline p}
$$
$$
= c_{p, q}\,\left(\int\limits_{0}^\infty f^{\underline p}(t)\,t^{\underline p- 1}\,\left\|\chi_{(0, \,x)}(t)\, \left[\frac{\omega_2(x)}{x}\right]^{\underline p}\right\|_{L_{\frac{q(\cdot)}{\underline p}}(0, \infty)}\,dt\right)^{1/\underline p}
$$
$$
= c_{p, q}\,\left(\int\limits_{0}^\infty f^{\underline p}(t)\,t^{\underline p- 1}\,\left\|\frac{\omega_2}{x}\right\|_{L_{q(\cdot)}(t, \infty)}^{\underline p}\,dt\right)^{1/\underline p}= c_{p, q}\,\left\|f\,t^{ 1/{\overline p'}}\,\left\|\frac{\omega_2}{x}\right\|_{L_{q(\cdot)}(t, \infty)} \right\|_{L_{\underline p}(0, \infty)}.
$$
Finally, apply Theorem 2, we get
$$
\left\|f\,t^{ 1/{\overline p'}}\,\left\|\frac{\omega_2}{x}\right\|_{L_{q(\cdot)}(t, \infty)} \right\|_{L_{\underline p}(0, \infty)}\le d_p\, \left\|\frac{t^{1/\overline p'}\, \left\|\frac{\omega_2}{x} \right\|_{L_{q(\cdot)}(t, \infty)}}{\omega_1} \right\|_{L_{r(\cdot)}(0,\infty)}\,\|f\|_{L_{p(\cdot),\,\omega_1}(0, \infty)}.
$$
Thus
$$
\left\|Hf\right\|_{L_{q(\cdot),\,\omega_2}(0, \infty)}\le {\underline p}^{\frac 1{\underline p}}\;c_{p, q}\;
d_p\, \left\|\frac{t^{1/p'}\,\left\|\frac{\omega_2}{x}\right\|_{L_{q(\cdot)}(t, \infty)}}{\omega_1}\right\|_{L_r(\cdot)(0, \infty)}\,\|f\|_{L_{p(\cdot),\,\omega_1}(0, \infty)}.
$$

The Theorem 6 is proved.
\begin{theorem} Let $0< \underline p\le p(x)\le q(x)\le  \overline q< 1,$ $\displaystyle{r(x)= \frac{\underline p\,p(x)}{p(x)- \underline p}}$ and $f(x)$ are non-negative and increasing function defined on $(0, \infty).$
Suppose $\omega_1$ and $\omega_2$ are weight functions defined on $(0, \infty).$

Then for any $f\in L_{p(x),\,\omega_1}(0,\,\infty)$ the inequality
$$
\left\|Hf\right\|_{L_{q(\cdot),\,\omega_2}(0, \infty)}\le {\underline p}^{\frac 1{\underline p}}\;c_{p,q}\,d_p\, \left\|\left\|\frac{(x- t)^{ 1/\overline p'}\,\omega_2}{x} \right\|_{L_{q(\cdot)}(t, \infty)}\,\frac 1{\omega_1} \right\|_{L_{r(\cdot)}(0,\infty)}\,\left\|f\right\|_{L_{p(\cdot),\,\omega_1}(0, \infty)},
$$
where $c_{p,q}$ and $d_p$ the constants in Theorem 6.
\end{theorem}

{\bf Proof.} Taking $a= 0,$ $b= x$ and $s= \underline p$ and apply Lemma 3 and property 5), we have
$$
\left\|Hf\right\|_{L_{q(\cdot),\,\omega_2}(0, \infty)}= \left\|\omega_2 Hf\right\|_{L_{q(\cdot)}(0, \infty)}= \left\|\frac{\omega_2}{x} \int\limits_{0}^x f(t)\,dt\right\|_{L_{q(\cdot)}(0, \infty)}
$$
$$
\le \left(\underline p\right)^{1/\underline p}\left\|\frac{\omega_2(x)}{x} \left(\int\limits_{0}^x f^{\underline p}(t)\,(x- t)^{\underline p- 1}\,dt\right)^{1/\underline p}\right\|_{L_{q(\cdot)}(0, \infty)}.
$$
Now applied Theorem 3, we get
$$
\left\|\frac{\omega_2(x)}{x}\,\left(\int\limits_{0}^x f^{\underline p}(t)\,(x- t)^{\underline p- 1}\,dt\right)^{1/\underline p}\right\|_{L_{q(\cdot)}(0, \infty)}
$$
$$
=\left\|\left(\int\limits_{0}^\infty f^{\underline p}(t)\,\chi_{(0, \,x)}(t)\,\left[\frac{\omega_2(x)}{x}\right]^{\underline p}\,(x- t)^{\underline p- 1}\,dt\right)^{1/\underline p}\right\|_{L_{q(\cdot)}(0, \infty)}
$$
$$
= \left\|\int\limits_{0}^\infty f^{\underline p}(t)\,\chi_{(0,\,x)}(t)\,\left[\frac{\omega_2(x)}{x}\right]^
{\underline p}\,(x- t)^{\underline p- 1}\,dt\right\|_{L_{\frac{q(\cdot)}{\underline p}}(0, \infty)}^{1/\underline p}
$$
$$
\le c_p\,\left(\int\limits_{0}^\infty \left\| f^{\underline p}(t)\,\chi_{(0, \,x)}(t)\,\left[\frac{\omega_2(x)}{x}
\right]^{\underline p}\,(x- t)^{\underline p- 1}\right\|_{L_{\frac{q(\cdot)}{\underline p}}(0, \infty)}\,dt\right)^{1/\underline p}
$$
$$
= c_p\,\left(\int\limits_{0}^\infty f^{\underline p}(t)\,\left\|\chi_{(0, \,x)}(t)\, \left[\frac{(x- t)^{ 1/\overline p'}}{x}\,\omega_2(x)\right]^{\underline p}\right\|_{L_{\frac{q(\cdot)}{\underline p}}(0, \infty)}\,dt\right)^ {1/\underline p}
$$
$$
= c_p\,\left(\int\limits_{0}^\infty f^{\underline p}(t)\,\,\left\|\frac{(x- t)^{ 1/\overline p'}}{x}\,\omega_2\right\|_{L_{q(\cdot)}(t, \infty)}^{\underline p}\,dt\right)^{1/\underline p}
$$
$$
= c_p\,\left\|f\,\left\|\frac{(x- t)^{ 1/\overline p'}}{x} \,\omega_2
\right\|_{L_{q(\cdot)}(t, \infty)} \right\|_{L_{\underline p}(0, \infty)}.
$$
Finally, apply Theorem 2, we get
$$
\left\|f\,\left\|\frac{(x- t)^{ 1/\overline p'}}{x}\,\omega_2\right\|_{L_{q(\cdot)}(t, \infty)} \right\|_{L_{\underline p}(0, \infty)}
$$
$$
\le \, \left\|\left\|\frac{(x- t)^{ 1/\overline p'}\,\omega_2}{x} \right\|_{L_{q(\cdot)}(t, \infty)}\,\frac 1{\omega_1} \right\|_{L_{r(\cdot)}(0,\infty)}\,\|f\|_{L_{p(\cdot),\,\omega_1}(0, \infty)}.
$$
Thus
$$
\left\|Hf\right\|_{L_{q(\cdot),\,\omega_2}(0, \infty)}\le {\underline p}^{\frac 1{\underline p}}\;c_{p,q}\,d_p\, \left\|\left\|\frac{(x- t)^{ 1/\overline p'}\,\omega_2}{x} \right\|_{L_{q(\cdot)}(t, \infty)}\,\frac 1{\omega_1} \right\|_{L_{r(\cdot)}(0,\infty)}\,\|f\|_{L_{p(\cdot),\,\omega_1}(0, \infty)}.
$$

The Theorem 7 is proved.

For the dual operator $H^*$ a theorem below is proved analogously.

\begin{theorem} Let $x\in (0, \infty),$ $0< \underline p\le p(x)\le q(x)\le  \overline q< 1,$ $\displaystyle{r(x)= \frac{\underline p\,p(x)}{p(x)- \underline p}}$ and $f(x)$ are non-negative and decreasing function defined on $(0, \infty).$ Suppose $\omega_1$ and $\omega_2$ are weight functions defined on $(0, \infty).$

Then for any $f\in L_{p(x),\,\omega_1}(0,\,\infty)$ the inequality
$$
\left\|H^*f\right\|_{L_{q(\cdot),\,\omega_2}(0, \infty)}\le {\underline p}^{\frac 1{\underline p}}\,c_{p,q}\,d_p\, \left\|\left\|\frac{(t- x)^{ 1/\overline p'}\,\omega_2}{x} \right\|_{L_{q(\cdot)}(0,\,t)}\,\frac 1{\omega_1} \right\|_{L_{r(\cdot)}(0,\infty)}\,\left\|f\right\|_{L_{p(\cdot),\,\omega_1}(0, \infty)},
$$
where $c_{p,q}$ and $d_p$ the constants in Theorem 6.
\end{theorem}

\begin{remark} Note that Theorem 6,Theorem 7 and Theorem 8 in the case $p(x)= q(x)= p= const$ and
$\omega_1(x)= \omega_2(x)= x^\alpha$ was proved in [6] (see also [5]). In the case $1\le p(x)\le
q(x)\le \overline q< \infty$ Hardy inequality is very much studied (see [2], [3] and etc.).
In the constant exponent case $1\le p\le q\le \overline q\le \infty$ for detailed information we
refer to [10]. Note that similar problem for Hardy maximal function was investigated in [9] and [11].
\end{remark}
{\bf Example 4.2.} {\it Let $x\in (0,\,\infty),\,$ $0< p(x)= p= const< 1,$
$q(x)= \left\{
             \begin{array}{l}
            \frac 14, \;\;\, for\; 0< x< 1\\
            \frac 12, \;\;\, for \; x\ge 1,\\
            \end{array}
            \right.
$
$0< p\le q(x)$ and $\displaystyle{p'= \frac{p}{p- 1}}.$
Suppose $\omega_1(x)= x^{\alpha},$ $\omega_2(x)= x^{\beta+ 1},$ $\beta< -2,$ $\beta\ne -4$ and
$\displaystyle{\beta+ 2+ \frac 1{p'}< \alpha< \min\left\{\frac 1{p'};\;\beta+ 4+ \frac 1{p'}\right\} },$ where
$r(x)= \infty.$

Then the pair $\left(\omega_1, \omega_2\right)$ satisfies the condition of Theorem 6.}
\\
{\bf Example 4.3.} {\it Let $x\in (0,\,\infty),\,$ $0< \underline p\le p(x)\le q(x)\le \overline q< 1$ and $\displaystyle{\overline p'= \frac{\underline p}{\underline p- 1}}.$
Suppose $\displaystyle{\omega_1(x)= x^{1/\overline{p}'}\left\|\frac{\omega_2}{x}\right\|_{L_{q(\cdot)}(x, \infty)}}.$  Then condition $\displaystyle{\|1\|_{L_{r(\cdot)}(0, \infty)}< \infty}$ is guaranteed the satisfy of condition of Theorem 6.
Note that by Definition 1 the condition $\displaystyle{\|1\|_{L_{r(\cdot)}(0, \infty)}< \infty}$ is equivalent to
$$
\int\limits_0^\infty \delta^{\frac{\underline p\,p(x)}{p(x)- \underline p}}\,dx< \infty,
$$
where $\delta\in (0, 1).$
Then the pair $\left(\omega_1, \omega_2\right)$ satisfies the condition of Theorem 6.}

{\bf Acknowledgement.} This work was supported by the Science Development Foundation under the President of
the Republic of Azerbaijan EIF-2010-1(1)-40/06-1.

\begin{center}
{\bf References}
\end{center}

[1] R.A.Bandaliev, {\it On an inequality in Lebesgue space with
mixed norm and with variable summability exponent,} Mat. Zametki, {\bf 3}
(84)(2008), 323-333.(In Russian). English translation:  Math. Notes,
{\bf 3}(84)(2008), 303-313 (2008).

[2] R.A.Bandaliev, {\it The boundedness of certain sublinear
operator in the weighted variable Lebesgue spaces,}  Czechoslovak Math. J.
{\bf 60}(2), 327-337 (2010).

[3] R.A.Bandaliev, {\it The boundedness of multidimensional Hardy operator in the weighted variable
Lebesgue spaces}, Lithuanian Math. J. {\bf 50}(2010), no.3, 249-259.

[4] R.A.Bandaliev, Z.V.Safarov, {\it Criteria of two-weighted inequalities for multidimensional
Hardy type operators in weighted Musielak-Orlicz spaces and some applications,} Mathematische
Nachrichten, 2012 (accepted).

[5] R.A.Bandaliev, {\it Embedding between variable exponent Lebesgue spaces with measures,}
Azerbaijan Journal of Math., {\bf 2}(1)(2012), 111-117.

[6] R.A.Bandaliev and K. K. Omarova, {\it Two-weight norm inequalities for certain singular integrals,} Taiwanese
Journal of Math.,{\bf 2} (2012), 113-132.

[7] J.Bergh, V.I.Burenkov, L.-E. Persson, {\it On some sharp reversed H\"{o}lder
and Hardy-type inequalities,} Math. Nachr., 169 (1994), 19-29.

[8] V.I. Burenkov, {\it On the exact constant in the Hardy inequality with 0 < p < 1 for monotone
functions,} Trudy Matem. Inst. Steklov. 194 (1992), 58-62 (in Russian). English transl. in Proc.
Steklov Inst. Math., 194, no. 4 (1993), 59-63.

[9] L.Diening, P.Harjulehto, P.H\"{a}st\"{o},   and  M. R\.{u}\v{z}i\v{c}ka,
{\it Lebesgue and Sobolev spaces with variable exponents,} Springer Lecture Notes, v.2017,
Springer-Verlag, Berlin, 2011.

[10] O.Kov\'{a}\v{c}ik, J. R\'{a}kosn\'{\i}k, {\it On spaces $L^{p(x)}$ and $W^{k, p(x)},$}
Czechoslovak Math. J. ({\bf 41}){\bf 116} (1991) 592-618.

[11] A.K.Lerner, {\it On some questions related to the maximal operator on variable $L^p$ spaces,}
Trans. Amer. Math. Soc., {\bf 362}(2010), no. 8, 4229-4242.

[12] V.G.Maz'ya, {\it Sobolev spaces,} (Springer-Verlag, Berlin, 1985).

[13] B.Muckenhoupt, {\it Weighted norm inequalities for the Hardy maximal function,} Trans. Amer. Math. Soc.,
{\bf 166}(1972).

[14] J. Musielak, {\it Orlicz spaces and modular spaces,} Lecture Notes in Math.1034. Springer-Verlag,
Berlin-Heidelberg-New York, 1983.

[15] W. Orlicz, {\it \"{U}ber konjugierte exponentenfolgen,}  Studia Math.{\it 3}(1931) 200-212.

[16] K.R.Rajagopal, M. R\.{u}\v{z}i\v{c}ka, {\it Mathematical modeling of  electrorheological
materials,} Cont. Mech. and Termodyn., {\bf 13}(2001) 59-78.

[17] S.G.Samko. {\it "Differentiation and integration of variable order and the spaces $L^{p(x)}$",}
\linebreak {\it Proc.Inter.Conf "Operator theory for complex and hypercomplex analysis",}
Mexico, 1994,  {\it Contemp. Math.}, {\bf 212}(1998), 203-219.

[18] I.I.Sharapudinov, {\it On a topology of the space $L^{p(t)}([0,1]),$}
{\it Matem. Zametki,} {\bf 26,} 613-632 (1979) (in Russian): English
translation: {\it Math. Notes,} {\bf 26}, 796-806 (1979).

[19] Q.H.Zhang, {\it Existence and asymptotic behavior of positive solutions for variable
exponent elliptic systems,} Nonlinear Analysis TMA, (1){\bf 70}(2009) 305-316.

[20] V.V. Zhikov, {\it Averaging of functionals of the calculus of variations and elasticity
theory,} Izv.  Akad. Nauk SSSR.{\bf 50}(1986) 675-710. (In Russian). English transl.:
Math. USSR, Izv., {\bf 29}(1987) 33-66.

\begin{flushleft}
DEPARTMENT OF MATHEMATICAL ANALYSIS, INSTITUTE OF MATHEMATICS \\
AND MECHANICS OF NATIONAL ACADEMY OF SCIENCES OF AZERBAIJAN,\\
Baku, Az 1141, B.Vahabzade str., 9\\
{\it E-mail address}: bandaliyev.rovshan@math.ab.az
\end{flushleft}

\end{document}